\documentclass[leqno]{amsart}
\usepackage{amsmath,amssymb,latexsym}
 \newtheorem{definition}{Definition}[section]
 \newtheorem{theorem}[definition]{Theorem}
 \newtheorem{lemma}[definition]{Lemma}
 \newtheorem{proposition}[definition]{Proposition}

 \newtheorem*{theorem*}{Theorem}
\newtheorem*{proposition*}{Proposition}
\newtheorem*{lemma*}{Lemma}

 \theoremstyle{remark}
 
 \newtheorem{remark}[definition]{Remark}

\newcommand{\op}[1]{\operatorname{#1}}

\newcommand{\acou}[2]{\ensuremath{\langle #1 , #2 \rangle}} 

\newcommand{\tr}{\ensuremath{\op{tr}}}
\newcommand{\Tr}{\ensuremath{\op{Tr}}}

\newcommand{\End}{\ensuremath{\op{End}}}

\newcommand{\Cl}{\ensuremath{{\op{Cl}}}}

\newcommand{\C}{\ensuremath{\mathbb{C}}} 
 
\newcommand{\R}{\ensuremath{\mathbb{R}}} 
\newcommand{\Z}{\ensuremath{\mathbb{Z}}}

\newcommand{\UR}{U\times\R}

\newcommand{\Ca}[1]{\ensuremath{\mathcal{#1}}}

\newcommand{\cE}{\Ca{E}}

\newcommand{\cS}{\ensuremath{\mathcal{S}}}
\newcommand{\cD}{\ensuremath{\mathcal{D}}}

\newcommand{\pdo}{$\Psi$DO}
\newcommand{\pvdo}{\ensuremath{\Psi_{\op{v}}}} 
\newcommand{\psido}{$\Psi$DO} 
\newcommand{\psidos}{$\Psi$DO's}

\newcommand{\ev}{{\text{ev}}}
\newcommand{\odd}{{\text{odd}}}

\newcommand{\hotimes}{\hat\otimes}

\newcommand{\sD}{\ensuremath{{/\!\!\!\!D}}}

\newcommand{\sQ}{\ensuremath{{/\!\!\!\!Q}}}

\newcommand{\sS}{\ensuremath{{/\!\!\!\!\!\;S}}}

\newcommand{\clD}{\ensuremath{\sD}}
\newcommand{\clQ}{\ensuremath{Q}}

\newcommand{\clR}{\ensuremath{R}}

\newcommand{\sDp}{\ensuremath{\sD_{\C^{p}}}}
\newcommand{\sQp}{\ensuremath{\tilde{Q}}}

\begin{document}
\title{A new proof of the local regularity of the eta invariant of a Dirac operator} 

\author{Rapha\"el Ponge}

\address{Department of Mathematics, Ohio State University, Columbus, USA.}
\email{ponge@math.ohio-state.edu}
\thanks{Research partially supported by NSF grant DMS 0409005}

\numberwithin{equation}{section}

\date{}

\begin{abstract}
In this paper we use the approach of~\cite{Po:CMP1} to the proof of the Local Index Theorem to give a new proof of Bismut-Freed's 
result on the local regularity of the eta invariant of a Dirac operator in odd dimension.   
\end{abstract}

\maketitle 

\section{Introduction}

The eta invariant of a selfadjoint elliptic \psido\ was introduced by 
    Atiyah-Patodi-Singer~\cite{APS:SARG1} as a boundary correction to their index formula on manifolds with boundary. It is obtained as the regular 
    value at $s=0$ of the eta function,  
    \begin{equation}
         \eta(P;s)=\Tr P|P|^{-(s+1)}=\int_{M} \eta(P;s)(x). 
    \end{equation}
    However, the residue at $s=0$ of the local eta function $\eta(P;s)(x)$ needs not vanish (see, e.g.,~\cite{Gi:RLEFO}), 
    so it is a nontrivial fact that the regular value exists~(\cite{APS:SARG3}, 
    \cite{Gi:RGEFO}). Nevertheless, global $K$-theoretic arguments  allows us to reduce the proof to the case of a Dirac operator on an odd dimensional 
    spin Riemannian manifold with coefficients in a Hermitian vector bundle, for which the result can be obtained by using invariant theory
    (see~\cite{APS:SARG3}, \cite{Gi:ITHEASIT}). 
    
      Subsequently, Wodzicki~(\cite{Wo:SAZF}, \cite{Wo:LISA}) generalized the result of Atiyah-Patodi-Singer and Gilkey to the nonselfadjoint 
      setting. More precisely, he proved that:\smallskip 
      
      (i) The regular value at $s=0$ of the zeta function $\zeta_{\theta}(P;s)=\Tr P_{\theta}^{-s}$ of an elliptic \psido\ is independent of the spectral cutting 
      $\{\arg \lambda =\theta\}$ used to define $P_{\theta}^{s}$ (see~\cite{Wo:SAZF}, \cite[1.24]{Wo:LISA});\smallskip
      
   (ii) The noncommutative residue of a \psido\ projection is always zero (see~\cite[7.12]{Wo:LISA}).\smallskip
   
   \noindent The original proofs of Wodzicki 
   are quite involved, but it follows from an observation of Br\"uning-Lesch~\cite[Lem.~2.6]{BL:OEICNLBVP} (see also~\cite[Rmk.~4.5 ]{Po:JFA1}) 
   that Wodzicki's results can be deduced from the aforementioned result of Atiyah-Patodi-Singer and Gilkey. 
    
    Next, in the case of a Dirac  operator $\sD_{\cE}$ on a odd spin Riemannian manifold $M^{n}$ with coefficients in a Hermitian bundle $\cE$, 
    Bismut-Freed~\cite{BF:AEF2} proved in a purely analytic fashion that the local eta function $\eta(P;s)(x)$ is actually regular at $s=0$. 
    Thanks to the Mellin formula for $\Re s >-1$ we have  
    \begin{equation}
        \sD_{\cE}|\sD_{\cE}|^{-(s+1)}=\Gamma(\frac{s+1}{2})^{-1}\int_{0}^{\infty}t^{\frac{s-1}{2}}\sD_{\cE}e^{-t\sD_{\cE}^{2}}dt.
%         \qquad .
%         \label{eq:}
    \end{equation}
    
    For $t>0$ let $h_{t}(x,y) \in C^{\infty}(M,|\Lambda|(M)\otimes \End \cE)$ be the 
    kernel of $\sD_{\cE}e^{-t\sD_{\cE}^{2}}$, where $|\Lambda|(M)$ denotes the bundle of densities on $M$. Since by standard heat kernel asymptotics we have
    $\tr h_{t}(x,x)=\op{O}(t^{-\frac{n}{2}})$ as $t\rightarrow 0^{+}$ (see Theorem~\ref{thm:Greiner.heat-kernel-asymptotics} ahead),  for $\Re s >n-1$ we get
%      for $ $ we get  
    \begin{equation}
        \eta(\sD_{\cE};s)(x)= \Gamma(\frac{s+1}{2})^{-1}\int_{0}^{\infty}t^{\frac{s-1}{2}}\tr h_{t}(x,x)dt.
    \end{equation}
   Then Bismut-Freed~\cite[Thm.~2.4]{BF:AEF2} proved that in the $C^{0}$-topology we have
    \begin{equation}
        \tr h_{t}(x,x) =\op{O}(\sqrt{t}) \qquad \text{as $t\rightarrow 0^{+}$}. 
         \label{eq:Intro.BFA}
    \end{equation}
    It thus follows that the local eta function $\eta(\sD_{\cE};s)(x)$ is actually holomorphic for $\Re s>-2$. In particular, we have the formula, 
    \begin{equation}
        \eta(\sD_{\cE})=\frac{1}{\sqrt{\pi}}\int_{0}^{\infty}t^{\frac{-1}{2}}\Tr \sD_{\cE}e^{-t\sD_{\cE}^{2}}dt,
    \end{equation}
    which, for instance, plays a crucial role in the study of the adiabatic limit of the eta invariant of a Dirac operator (see~\cite{BC:EIAL}, 
    \cite{BF:AEF2}). Incidentally, Bismut-Freed's asymptotics~(\ref{eq:Intro.BFA}), which is a purely analytic statement, implies the global regularity of 
    the eta invariant of any selfadjoint elliptic \psido, as well as the aforementioned generalizations of Wodzicki. 
    
   Now, the standard proof of the asymptotics~(\ref{eq:Intro.BFA}) is essentially based on a reduction to the local index theorem of Patodi, Gilkey, 
   Atiyah-Bott-Patodi~(\cite{ABP:OHEIT},  \cite{Gi:ITHEASIT}; see also~\cite{Ge:SPLASIT}), which  provides us with a heat kernel proof of the
   index theorem of Atiyah-Singer~(\cite{AS:IEO1}, \cite{AS:IEO3}) for Dirac operators. In the original proof of~(\ref{eq:Intro.BFA}) in~\cite{BF:AEF2} 
   the reduction is done by 
   introducing an extra Grassmanian variable $z$, $z^{2}=0$, and in~\cite[Sect.~8.3]{Me:APSIT} by means of a suspension argument. Moreover, 
   in~\cite{BF:AEF2} the authors refer to the results of~\cite{Gr:AEHE} to justify the differentiability of the heat kernel asymptotics.

   On the other hand, in~\cite{Po:CMP1} the approach to the heat kernel asymptotics of~\cite{Gr:AEHE} was combined with general 
   considerations on Getzler's order of Volterra \psidos\ to produce a new short proof of the local index theorem which holds in many other settings. 
   Furthermore, the arguments used in this proof have other applications such as the computation of the CM cyclic cocycle of~\cite{CM:LIFNCG} for Dirac operators. 
  
  The aim of this paper is to show that we can get a direct proof of Bismut-Freed's asymptotics by using the approach of~\cite{Po:CMP1}.
 In fact, once the background from~\cite{Gr:AEHE} and~\cite{Po:CMP1} is set-up, the proof 
  becomes extremely simple since it is deduced by combining the considerations on Getzler orders of~\cite{Po:CMP1} 
   with the observation that the subleading terms (in the Getzler order sense) 
   of the various asymptotics at stake vanish at first order. 

   As in~\cite{Po:CMP1} the approach of this paper is rather general and should therefore hold in various other settings. Moreover, it is believed that
   this approach could also be useful to study 
   adiabatic-type limits of eta invariants. For instance, it has been shown by Rumin~\cite{Ru:SRLDFSCM} that, under the so-called subriemannian limit, the 
   differential form spectrum of the de Rham complex on a contact manifold converges to that of the contact complex~\cite{Ru:FDVC}, 
   so it could be fruitful to use the approach of this paper to further study the 
   behavior under the subriemannian limit of the eta invariant on a contact manifold. 
    
    This paper is organized as follows. In Section~\ref{sec.volterra} we recall Greiner's approach to the heat kernel asymptotics.  In 
    Section~\ref{sec.getzler}, after having recalled the background of~\cite{Po:CMP1}, we prove Bismut-Freed's asymptotics.  

\section{Greiner's approach of the heat kernel asymptotics}
\label{sec.volterra}
  In this section we recall Greiner's approach of the heat kernel asymptotics as in \cite{Gr:AEHE} and 
\cite{BGS:HECRM}.% (see also~\cite{Po:CMP1}). 

Let $M^n$  be a compact Riemannian manifold, let   $\cE$  a Hermitian 
vector bundle over $M$ and  let $\Delta:C^{\infty}(M,\cE)\rightarrow C^{\infty}(M,\cE)$ be a selfadjoint second order elliptic differential operator
with principal symbol $a_{2}(x,\xi)>0$. Then $\Delta$ is bounded  from below on $L^{2}(M,\cE)$ and by standard functional calculus 
we can define 
 $e^{-t\Delta}$, $t\geq0$, as a selfadjoint bounded 
operator on $L^{2}(M,\cE)$. 
In fact, $e^{-t\Delta}$ is smoothing for 
$t>0$, so its Schwartz  kernel $k_{t}(x,y)$ is in  
$C^\infty(M,\cE)\hotimes C^\infty(M,\cE^{*}\otimes |\Lambda|(M))$,  
where $|\Lambda|(M)$ denotes the bundle of densities on $M$. 

Recall that the heat semigroup allows us to invert the heat equation, in the sense that the operator, 
\begin{equation}
    Q_{0}u(x,t)=\int_{0}^\infty e^{-s\Delta} u(x,t-s)ds, \qquad u \in C^\infty_{c}(M\times \R, \cE), 
     \label{eq:volterra.inverse-heat-operator}
\end{equation}
maps continuously into $C^{0}(\R, L^{2}(M,\cE)) \subset \cD'(M\times\R, \cE)$ and satisfies
\begin{equation}
    (\Delta+\partial_{t})Q_{0}u = Q_{0}(\Delta+\partial_{t})u=u \qquad \forall u \in C^\infty_{c}(M\times\R,\cE).
\end{equation}

Notice that the operator $Q_{0}$  has the \emph{Volterra property} in the sense of~\cite{Pi:COPDTV}: it is translation invariant and satisfies 
the causality principle, i.e., $Q$ has a 
distribution kernel of the form $K_{Q_{0}}(x,y,t-s)$ where $K_{Q_{0}}(x,y,t)$ vanishes on the region $t<0$. In fact, we have
\begin{equation}
    K_{Q_{0}}(x,y,t) = \left\{ 
    \begin{array}{ll}
         k_{t}(x,y) & \quad \text{if $t> 0$},  \\
        0 &  \quad \text{if $t<0$}. 
    \end{array}\right. 
\end{equation}

The above equalities lead us to use pseudodifferential techniques to study the heat kernel $k_{t}(x,y)$. 
The idea, which goes back to Hadamard, is to consider a class of \psido's, the Volterra \psido's 
(\cite{Gr:AEHE},~\cite{Pi:COPDTV}, \cite{BGS:HECRM}), 
taking into account:  \smallskip 

(i)  The aforementioned Volterra property;\smallskip

(ii) The parabolic homogeneity of the heat operator $\Delta+ \partial_{t}$, i.e.,  the homogeneity with respect 
             to the dilations, 
             \begin{equation}
                 \lambda.(\xi,\tau)=(\lambda\xi,\lambda^{2}\tau), \qquad (\xi,\tau)\in \R^{n+1}, \quad \lambda\neq 0. 
                 \label{eq:Greiner.parabolic-dilation}
             \end{equation}
             
    In the sequel for $g\in \cS'(\R^{n+1})$  and $\lambda\neq 0$ we let $g_{\lambda}$ be the tempered distribution 
    defined by   
    \begin{equation}
        \acou{g_{\lambda}( \xi,\tau)}{u(\xi,\tau)} =    |\lambda|^{-(n+2)} 
            \acou{g(\xi,\tau)} {u(\lambda^{-1}\xi, \lambda^{-2}\tau)}, \quad u \in \cS(\R^{n+1}). 
            \label{eq:Greiner.parabolic-dilations-distributions}
    \end{equation}
We then say that $g$ is \emph{parabolic 
homogeneous} of degree $m$, $m\in \Z$, when  we have 
$g_{\lambda}=\lambda^m g$ for any $\lambda \neq 0$.  
    
%     \begin{definition}
%     A distribution $ g\in \cS'(\R^{n+1})$ 
% \end{definition}
Let $\C_{-}$ denote the complex halfplane $\{\Im \tau >0\}$ with closure $\overline{\C_{-}}\subset \C$. Then:  
\begin{lemma}[{\cite[Prop.~1.9]{BGS:HECRM}}] \label{lem:volterra.volterra-extension}
Let $q(\xi,\tau)\in C^\infty((\R^{n}\times\R)\setminus0)$ be a parabolic homogeneous 
symbol of degree $m$ such that: \smallskip 

(i) $q$ extends to a continuous function on $(\R^{n}\times\overline{\C_{-}})\setminus0$ in 
    such way to be holomorphic in the last variable when the latter is restricted to $\C_{-}$. \smallskip 

\noindent Then  there is a unique $g\in \cS'(\R^{n+1})$ agreeing with $q$ 
on $\R^{n+1}\setminus 0$ so that: \smallskip 

(ii) $g$ is parabolic homogeneous of degree $m$; \smallskip 

(iii) The inverse Fourier transform $\check g(x,t)$ vanishes for $t<0$. 
\end{lemma}

Let $U$ be an open subset of $\R^{n}$. We define Volterra \psidos\  on 
$U\times\R$ as follows. 

\begin{definition}
    $S_{\op v}^m(U\times\R^{n+1})$, $m\in\Z$,  consists of smooth functions $q(x,\xi,\tau)$ on 
    $U\times\R^n\times\R$ with an asymptotic expansion  $q \sim \sum_{j\geq 0} q_{m-j}$, where: \smallskip 
    
    - The symbol $q_{l}(x,\xi,\tau)\in C^{\infty}(U\times[(\R^n\times\R)\setminus0])$ is  parabolic homogeneous of degree $l$ and  satisfies the 
    property (i) in Lemma~\ref{lem:volterra.volterra-extension} with respect to the variables $\xi$ and $\tau$; \smallskip 
    
    - The sign $\sim$ means that, for any integer $N$ and any compact $K\subset U$, there is a constant $C_{NK\alpha\beta 
    k}>0$ such that, for  $x\in K$ and $|\xi|+|\tau|^{\frac12}>1$, we have 
            \begin{equation}
                |\partial^{\alpha}_{x}\partial^{\beta}_{\xi} \partial^k_{\tau}(q-\sum_{j< N} 
            q_{m-j})(x,\xi,\tau) | 
                \leq C_{NK\alpha\beta k} (|\xi|+|\tau|^{1/2})^{m-N-|\beta|-2k}. 
                          \label{eq:volterra.asymptotic-symbols}
            \end{equation}
\end{definition}

Given a symbol $q \in S_{\op v}^m(U\times\R^{n+1})$ we can quantize it as a standard symbol by associating to it the operator 
$Q(x,D_{x},D_{t})$ from $C_{c}^\infty(U_{x}\times\R_{t})$ to 
    $C^\infty(U_{x}\times\R_{t})$ such that, for any $u \in C_{c}^\infty(U_{x}\times\R_{t})$, we have
\begin{equation}
  Q(x,D_{x},D_{t})u(x,t)=(2\pi)^{-(n+1)}\int e^{i(x.\xi+t.\tau)}q(x,\xi,\tau)\check{u}(\xi,\tau)d\xi d\tau. 
%     \label{eq:¥}
\end{equation}
This is a continuous operator from $C_{c}^\infty(U_{x}\times\R_{t})$ to 
    $C^\infty(U_{x}\times\R_{t})$ satisfying the Volterra property and with distribution kernel $\check{q}_{(\xi,\tau)\rightarrow 
    (y,t)}(x,y,t-s)$.  

\begin{definition}\label{def:volterra.PsiDO}
    $\pvdo^m(U\times\R)$, $m\in\Z$,  consists of continuous operators 
    $Q$ from $C_{c}^\infty(U_{x}\times\R_{t})$ to 
    $C^\infty(U_{x}\times\R_{t})$ such that: \smallskip 
    
    (i) $Q$ has the Volterra property; \smallskip 
    
    (ii) $Q$ is of the form $Q=q(x,D_{x},D_{t})+R$ for some symbol $q$ in $S^m_{\op v}(U\times\R)$ and some smoothing operator  $R$. 
\end{definition}

In the sequel if $Q$ is a Volterra \psido\ we let $K_{Q}(x,y,t-s)$ denote its distribution kernel, so that the 
distribution $K_{Q}(x,y,t)$ vanishes for $t<0$.  

Examples of Volterra \psido\ includes differential operators, as well as homogeneous operators below.

% \begin{example}
%     Let $P$ be a differential operator of order $2$ on $U$ and let $p_{2}(x,\xi)$ denote the principal symbol of $P$. Then 
%     the heat operator $P+\partial_{t}$ is a Volterra \psido\ of order $2$ with principal symbol $p_{2}(x,\xi)+i\tau$. 
% \end{example}
% 
% Other examples of Volterra \psido's are given by the 
\begin{definition}\label{def:Greiner.homogeneous-PsiDO} 
Let $q_{m}(x,\xi,\tau) \in C^\infty(U\times(\R^{n+1}\setminus 0))$ be a homogeneous Volterra 
symbol of order $m$ and let $g_{m}\in C^\infty(U)\hotimes 
    \cS'(\R^{n+1})$ denote its unique homogeneous extension given by 
    Lemma~\ref{lem:volterra.volterra-extension}. Then: \smallskip 
    
    -  $\check q_{m}(x,y,t)$ denotes the inverse Fourier transform $(g_{m})^{\vee}_{(\xi,\tau)\rightarrow (y,t)}(x,y,t)$;\smallskip 
% 	in the last $n+1$ variables;  
    
        -  $q_{m}(x,D_{x},D_{t})$ denotes the operator with kernel $\check q_{m}(x,y-x,s-t)$.
\end{definition}
\begin{remark}
    The above definition makes sense since it follows from the proof of Lemma~\ref{lem:volterra.volterra-extension} 
    in \cite{BGS:HECRM} that the extension process of 
    Lemma~\ref{lem:volterra.volterra-extension}   
    applied to each symbol $q_{m}(x,.,.)$, $x \in U$, is smooth with respect to $x$, so really gives rise to an element of $C^\infty(U)\hotimes 
    \cS'(\R^{n+1})$. 
 \end{remark}

\begin{proposition}[\cite{Gr:AEHE},~\cite{Pi:COPDTV}, \cite{BGS:HECRM}]\label{prop:Volterra-properties}
  The  following properties hold.\smallskip 

 \emph{1) Composition}. Let $Q_{j}\in \pvdo^{m_{j}}(\UR)$, $j=1,2$, have symbol $q_{j}$ 
and 
    suppose that $Q_{1}$ or $Q_{2}$ is properly supported. Then $Q_{1}Q_{2}$ belongs  to $\pvdo^{m_{1}+m_{2}}(\UR)$
    and has symbol  $q_{1}\#q_{2} \sim \sum \frac{1}{\alpha!} 
\partial_{\xi}^{\alpha}q_{1} D_{x}^\alpha q_{2}$.\smallskip 

\emph{2) Parametrices.} An operator $Q\in \pvdo^{m}(U\times\R)$ has a  parametrix 
    in $\pvdo^{-m}(U\times\R)$ if, and only if, its principal symbol is nowhere vanishing on 
    $U\times[(\R^n\times \overline{\C_{-}}\setminus 0)]$.\smallskip
    
    \emph{3) Invariance.} Let $\phi : U \rightarrow V$ be 
    a diffeomorphism onto another open 
subset $V$ of $\R^n$ and let $Q$ be a Volterra \pdo\ on $U\times\R$ of order $m$. Then 
$\clQ=(\phi\oplus \op{id}_{\R})_{*}Q$  is a Volterra \pdo\ on $V\times \R$ of order $m$. 
\end{proposition}
 \begin{remark}
 The proofs of the first and third properties above follow along the same lines as that of the corresponding proofs for standard \psidos. 
 For the second one we need the Volterra calculus to be asymptotically complete. The latter 
 fact cannot, however, be reached by means of the standard proof for classical \psidos, because the cut-off arguments therein are not valid
  valid anymore with analytic symbols. Nevertheless, simple proofs can be found in~\cite{Po:JAM1}.
 \end{remark}

The key property of Volterra \psidos\ which allows us to derive a differentiable heat kernel asymptotics is the following.
\begin{lemma}[{\cite[Chap.~I]{Gr:AEHE}}]\label{lem:volterra.key-lemma}
Let $Q\in \pvdo^{m}(U\times\R)$ have symbol $q \sim \sum_{j \geq 0} q_{m-j}$. Then the following asymptotics 
holds in  $C^\infty(U)$, 
\begin{equation}
    K_{Q}(x,x,t) \sim_{t\rightarrow 0^{+}} t^{-(\frac{n}2+[\frac{m}2]+1)} \sum_{l\geq 0} t^l \check{q}_{2[\frac{m}2]-2l}(x,0,1), 
    \label{eq:volterra.asymptotics-Q}
\end{equation}
     where the notation $\check{q}_{k}$ has the same meaning as in Definition~\ref{def:Greiner.homogeneous-PsiDO}. 
\end{lemma}
\begin{proof}
   As the Fourier transform relates the decay at infinity to 
    the behavior at the origin of the Fourier transform the 
    distribution $\check{q} -\sum_{j\leq J} \check{q}_{m-j}$ is in $C^{N}
    (U_{x}\times\R^{n}_{y}\times\R_{t})$ 
    as soon as $J$ is large enough. Since $Q-q(x,D_{x},D_{t})$ is smoothing,  we see that 
   $R_{J}(x,t)=K_{Q}(x,x,t) -\sum_{j\leq J} \check q_{m-j}(x,0,t)$ is of class $C^{N}$. 
   As $R_{J}(x,y,t)=0$ for $t<0$ we get $\partial_{t}^l R_{J}(x,0)=0$ for 
    $l=0,1,\ldots,N$, so that $R_{J}(.,t)$ is a $\op{O}(t^{N})$ in $C^N(U)$ as $t\rightarrow 
    0^{+}$. It then follows that in $C^\infty(U)$ we have  
\begin{equation}
            K_{Q}(x,x,t) \sim \sum_{j\geq 0} t^{\frac{-n}2-1}\check{q}_{m-j}(x,0,t) \qquad \text{as $t\rightarrow 
            0^{+}$}.
\end{equation}
        
    Now, using~(\ref{eq:Greiner.parabolic-dilations-distributions}) we see that for $\lambda\neq 0$ we have
    \begin{equation}
        (\check{q}_{m-j})_{\lambda}= |\lambda|^{-(n+2)} (q_{m-j,\lambda^{-1}})^{\vee}=  |\lambda|^{-(n+2)} \lambda^{j-m} 
    \check{q}_{m-j}.
         \label{eq:Greiner.homogeneity-inverse-Fourier-transform}
    \end{equation}
   Setting $\lambda=\sqrt{t}$ with $t>0$ then gives
 $\check  q_{m-j}(x,0,t)=t^{\frac{j-n-m}2+1} \check q_{m-j}(x,0,1)$. Thus, 
\begin{equation}
    K_{Q}(x,x,t)\sim  \sum_{j\geq 0} t^{\frac{j-n-m}2-1} \check{q}_{m-j}(x,0,1) \qquad \text{as $t\rightarrow 
            0^{+}$}.
     \label{eq:volterra.asymptotics-Q2}
\end{equation}

On the other hand, setting $\lambda=-1$ in~(\ref{eq:Greiner.homogeneity-inverse-Fourier-transform}) 
shows that when $m-j$ is odd we have $\check{q}_{m-j}(x,0,1)=-\check{q}_{m-j}(x,0,1)=0$ . Therefore, 
in~(\ref{eq:volterra.asymptotics-Q2})  only the 
symbols of even degree contributes to the asymptotics, i.e., we get~(\ref{eq:volterra.asymptotics-Q}).
% \begin{equation}
%     K_{Q}(x,x,t) \sim t^{-(\frac{n}2+[\frac{m}2]+1)} \sum_{l\geq 
% 0} t^l \check{q}_{2[\frac{m}2]-2l}(x,0,1) \qquad \text{as $t\rightarrow 0^{+}$}. 
% %     \label{eq:¥}
% \end{equation}
% The lemma is thus proved.
\end{proof}

The invariance property in Proposition~\ref{prop:Volterra-properties}
allows us to define Volterra \pdo's on $M\times \R$ acting on the 
sections of the vector bundle $\cE$. All the preceding  properties hold \emph{verbatim} in this context. 
In particular the heat operator $\Delta+\partial_{t}$  admits a  parametrix  $Q$ in $\pvdo^{-2}(M,\times\R,\cE)$. 

In fact, comparing the operator~(\ref{eq:volterra.inverse-heat-operator}) with any 
Volterra \psido-parametrix for $\Delta+\partial_{t}$ allows us to prove:  
\begin{theorem}[\cite{Gr:AEHE},~\cite{Pi:COPDTV}, {\cite[pp.~363-362]{BGS:HECRM}}]\label{thm:volterra.inverse} 
    The differential operator $\Delta+\partial_{t}$ has an inverse $Q_{0}\in \Psi^{-2}(M\times \R,\cE)$, so that 
    $Q_{0}(\Delta+\partial_{t})=(\Delta+\partial_{t})Q_{0}=1$. 
\end{theorem}

We can now get the heat kernel asymptotics in the differentiable version below.  

\begin{theorem}[{\cite[Thm.~1.6.1]{Gr:AEHE}}]\label{thm:Greiner.heat-kernel-asymptotics}
    Let $P:C^{\infty}(M,\cE)\rightarrow C^{\infty}(M,\cE)$ be a differential operator of order $m$ and let $h_{t}(x,y)$ 
    denote the distribution kernel of $Pe^{-t\Delta}$. Then in $C^{\infty}(M,|\Lambda|(M)\otimes \End \cE)$ as $t\rightarrow 0^{+}$ we have 
    \begin{equation}
        h_{t}(x,x)\sim t^{-[\frac{m}{2}]-\frac{n}{2}} \sum_{l \geq 0} t^{l} b_{l}(P,\Delta)(x), \quad 
        b_{l}(P,\Delta)(x)= \check{q}_{2[\frac{m}{2}]-2-2l}(x,0,1),
         \label{eq:Greiner.heat-kernel-asymptotics}
    \end{equation}
    where the 
equality on the right-hand side gives a formula for computing the densities $b_{l}(P,\Delta)(x)$'s in local trivializing coordinates 
using the symbol  $q\sim \sum_{j \geq 0} q_{m-2-j}$ of  $Q=P(\Delta+\partial_{t})^{-1}$.
\end{theorem}
\begin{proof}
    As $h_{t}(x,y)=P_{x}k_{t}(x,y)=P_{x}K_{(\Delta+\partial_{t})^{-1}}(x,y,t)=K_{P(\Delta+\partial_{t})^{-1}}(x,y,t)$ 
    the result follows by applying Lemma~\ref{lem:volterra.key-lemma} to $P(\Delta+\partial_{t})^{-1}$. 
\end{proof}

\section{Proof of Bismut-Freed's asymptotics}
\label{sec.getzler}
 In this section we shall show that implementing the rescaling of Getzler~\cite{Ge:SPLASIT} into the Volterra calculus allows us to get a new proof 
 of Bismut-Freed's asymptotics~(\ref{eq:Intro.BFA}). 

 Let $(M^n,g)$ be a compact Riemannian manifold of odd dimension and let $\Cl(M)$ denote its Clifford bundle, so that the fiber $\Cl_{x}(M)$ at 
 $x\in M$ is the complex algebra generated by $1$ and $T^{*}_{x}M$ with relations,
\begin{equation}
    \xi.\eta + \eta.\xi = -2\acou\xi\eta,  \qquad \xi,\eta\in T^{*}_{x}M.
\end{equation}

For $\xi \in \Lambda T^{*}_{\C}M$ and $l=0,\ldots,n$ we let $\xi^{(l)}$ denote component of $\xi$ in $\Lambda^{l}T^{*}_{\C}M$. Then 
the quantization map $c:\Lambda T^{*}_{\C}M \rightarrow \Cl(M)$ and its inverse  
$\sigma=c^{-1}$ are such that for  $\xi$ and $\eta$ in $\Lambda T^{*}_{\C}M$ we have   
\begin{equation}
    \sigma(c(\xi^{(i)})c(\eta^{(j)}))= \xi^{(i)}\wedge\eta^{(j)} \quad \bmod 
\Lambda^{i+j-2}T^{*}_{\C}M. 
    \label{eq:getzler.clifford-multiplication}
\end{equation}

Now, let $\cE$ be a Clifford bundle over $M$ equipped with a unitary Clifford connection $\nabla^{\cE}$ and let  $\sD_{\cE}$ be the associated 
 Dirac operator, %with coefficients in $\cE$, %given by the composition
\begin{equation}
  \sD_{\cE}:  C^\infty(M, \cE) \stackrel{\nabla^{\cE}}{\longrightarrow} C^\infty(M,T^{*}M \otimes\cE) 
    \stackrel{c}{\longrightarrow} C^\infty(M,\cE),
     \label{eq:odd-chern.Dirac-operator}
\end{equation}
where  $c$ denotes the Clifford action of $\Lambda T^{*}M$ on $\cE$. Recall that by the Lichnerowicz's formula~\cite[Thm.~3.52]{BGV:HKDO}
we have
 \begin{equation}
 \sD^{2}_{\cE}=
(\nabla^{\cE}_{i})^{*}\nabla^{\cE}_{i} + \mathcal{F}^{\cE/\sS}
 + \frac{\kappa^{M}}4,    
    \label{eq:BFA.Lichnerowicz}
\end{equation}
where $\kappa^{M}$ denotes the scalar curvature of $M$ and $\mathcal{F}^{\cE/\sS}$  is the twisted curvature
 $F^{\cE/\sS}$  acting by Clifford multiplication on $\cE$. 

For $t>0$ let $h_{t}(x,y)$ denote the kernel of $\sD_{\cE}e^{-t\sD_{\cE}^{2}}$. The aim of this section is to prove 
Bismut-Freed's asymptotics in the version below. 

\begin{theorem}[{\cite[Thm.~2.4]{BF:AEF2}}]
In $C^{\infty}(M,|\Lambda|(M)\otimes \End \cE)$ we have 
 \begin{equation}
     \tr_{\cE}h_{t}(x,x)=\op{O}(\sqrt{t}) \qquad \text{as $t\rightarrow 0^{+}$}.
      \label{eq:BFA.main-asymptotics}
 \end{equation}
\end{theorem}

First, as by Theorem~\ref{thm:Greiner.heat-kernel-asymptotics} we already have an asymptotics
in $C^\infty(M, |\Lambda|(M))$ for $ \tr_{\cE} h_{t}(x,x)$ as $t\rightarrow 0^{+}$ it is enough to prove~(\ref{eq:BFA.main-asymptotics}) 
at a point $x_{0}\in M$. In fact, proving~(\ref{eq:BFA.main-asymptotics})  at $x_{0}$ is a purely local issue since in local trivializing  coordinates near $x_{0}$ 
the coefficients of the asymptotics for $ \tr_{\cE} h_{t}(x,x)$ depend only the 
homogeneous components of the symbol of 
$\sD_{\cE}(\sD_{\cE}^{2}+\partial_{t})^{-1}$. In particular, if in local trivializing coordinates near $x_{0}$ we let 
$\sQ_{\cE}$ be a Volterra \psido\ parametrix for 
$\sD_{\cE}^{2}+\partial_{t}$  then we have
\begin{equation}
    h_{t}(x_{0},x_{0})=K_{\sD_{\cE}\sQ_{\cE}}(x_{0},x_{0},t) +\op{O}(t^{\infty}) \qquad \text{as $t\rightarrow 0^{+}$}.
     \label{eq:BFA.htKDQt}
\end{equation}

As usual it will be convenient to use normal coordinates centered at $x_{0}$ and trivializations of $TM$ and $\cE$ by means of synchronous orthogonal frames, 
assuming that the synchronous tangent frame $e_{1},\ldots, e_{n}$ is such that $e_{j}=\partial_{j}$ at $x=0$. This allows us to:\smallskip 

- Replace $\sD_{\cE}$ by a Dirac operator $\sDp$ on $\R^{n}$ acting on a trivial twisted bundle with fiber 
$\sS_{n}\otimes \C^{p}$, where $\sS_{n}$ denotes the spinor space of $\R^{n}$. 

- Have a metric $g$ and coefficients $\omega_{ikl}= \acou{\nabla^{LC}_{i}e_{k}}{e_{l}}$ of the 
Levi-Civita connection with behaviors near $x=0$ of the form  
\begin{equation}
     g_{ij}(x)=\delta_{ij}+\op{O}(|x|^{2}), \qquad  \omega_{ikl}(x)= -\frac12 
R_{ijkl}^{M}(0)x^j 
     +\op{O}(|x|^{2}),   
    \label{eq:AS-asymptotic-geometric-data}
\end{equation}
where $R_{ijkl}^{M}(0)=\acou{R^{M}(0)(\partial_{i},\partial_{j})\partial_{k}}{\partial_{l}}$. 

Second, let $\Cl(n)$ and $\Lambda(n)$ respectively denote the Clifford algebra and the (complexified) exterior algebra of $\R^{n}$. 
Since the dimension $n$ is odd  if we regard $c(dx^{i_{1}})\cdots c(dx^{i_{k}})$, $i_{1}<\ldots<i_{k}$, 
as an endomorphism of $\sS_{n}$ then, as pointed out in~\cite[p.~107]{BF:AEF2}, we have 
\begin{equation}
    \tr_{\sS_{n}} c(dx^{i_{1}})\cdots c(dx^{i_{k}}) = \left\{  \begin{array}{ll}
                                                              2^{[\frac{n}2]} & \text{if $k=0$,}   \\
                                                              0 & \text{if $0<k< n$,}   \\
          (-i)^{[\frac{n}2]+1} 2^{[\frac{n}2]}& \text{if $k=n$}. 
    \end{array} \right. 
    \label{eq:BFA.trace-clifford}
\end{equation}
% In other words, as soon as an odd number of Clifford variables are involved, the trace on $\sS_{n}$ behaves like the supertrace in even dimension. 

Next, the Dirac operator $\sD_{\C^{p}}$ is obtained by composing the action of $\Cl(n)$ on $\sS_{n}$ with the Clifford Dirac operator 
$\clD=c(dx^{i})\nabla_{{i}}$ with coefficients in $\Cl^{\odd}(n)$. Since $\Cl^{\ev}(n)$ acts on itself by 
multiplication, any Volterra \psido-parametrix $\clQ$ for $\clD^{2}+\partial_{t}$ has coefficients in $\Cl^{\ev}(n)$ up to a smoothing operator.  
By composing $\clD \clQ$  
with the Clifford action we get  the operator $\sDp \sQp$, where $\sQp$ is the 
Volterra \psido-parametrix of $\sDp^{2}+\partial_{t}$ corresponding to $\clQ$. Since 
$\clD\clQ$ has coefficients in $\Cl^{\odd}(n)$ up to a smoothing operator, using~(\ref{eq:BFA.htKDQt}) and~(\ref{eq:BFA.trace-clifford}) 
we see that as $t\rightarrow 0^{+}$ we have
\begin{multline}
    \tr_{\cE}h_{t}(0,0)= \tr_{\sS_{n}\otimes \C^{p}}K_{\sDp \sQp}(0,0,t) +\op{O}(t^{\infty})\\  = 
    (-i)^{[\frac{n}2]+1} 2^{[\frac{n}2]} (\sigma\otimes \tr_{\C^{p}})[ K_{\clD \clQ}(0,0,t)]^{(n)} +\op{O}(t^{\infty}).
     \label{eq:BFA.trace-ht-sigma}
\end{multline}
Therefore, we are reduced to study  the small time behavior of
$ (\sigma\otimes \tr_{\C^{p}})[ K_{\clD \sQ}(0,0,t)]^{(n)}$.

Now, the Getzler's rescaling~\cite{Ge:SPLASIT} aims to assign the degrees,
\begin{equation}
  \deg \partial_{j}=\deg c(dx^j)= 1, \quad \deg \partial_{t}=2, \quad  \deg x^j =-1 , 
%   \deg \partial_{j}=\frac12 \deg \partial_{t}=\deg c(dx^j)=- \deg x^j =1 , \qquad 
  \label{eq:AS-Getzler-order}
\end{equation}
while $ \deg B=0$ for any $B\in M_{p}(\C)$. 
As pointed out in~\cite{Po:CMP1} this allows us to define a filtration on Volterra \psido's with coefficients in $\Cl(n)$ as follows.  

Let $Q\in  \Psi_{\op{v}}^{*}(\R^{n}\times\R, \C^{p})\otimes \Cl(n)$ have symbol $q(x,\xi,\tau) \sim \sum_{k\leq m'} 
q_{k}(x,\xi,\tau)$. Then taking components in each subspace 
 $\Lambda^j (n)$ and then using Taylor expansions at $x=0$ gives formal expansions  
\begin{equation}
    \sigma[q(x,\xi,\tau)] \sim \sum_{j,k} \sigma[q_{k}(x,\xi,\tau)]^{(j)} \sim \sum_{j,k,\alpha} 
    \frac{x^{\alpha}}{\alpha!}  \sigma[\partial_{x}^{\alpha}q_{k}(0,\xi,\tau)]^{(j)},
    \label{eq:AS.Getzler-asymptotic}
\end{equation}
where the last asymptotic expansion is taken with respect to the filtration of 
$S_{\op v}^k(\R^{n}\times\R^{n+1}, \C^{p})\otimes \Lambda(n)$ by the subspaces 
$x^{\alpha}S_{\op v}^m(\R^{n}\times\R^{n+1}, \C^{p})\otimes \Lambda^{j}(n)$.  

Motivated by~(\ref{eq:AS-Getzler-order}) we shall say that the symbol $ \frac{x^{\alpha}}{\alpha!}  
\partial_{x}^{\alpha}\sigma[q_{k}(0,\xi,\tau)]^{(j)}$ is Getzler homogeneous of degree $k+j-|\alpha|$. Thus,  
we have an asymptotic expansion, 
% $\sigma[q(x,\xi,\tau)]$ as 
\begin{equation}
    \sigma[q(x,\xi,\tau)] \sim \sum_{j\geq 0} q_{(m-j)}(x,\xi,\tau), \qquad q_{(m)}\neq 0, 
    \label{eq:index.asymptotic-symbol}
\end{equation}
where the symbol $q_{(m-j)}$ is Getzler homogeneous of degree $m-j$. 

We shall call $m$ the \emph{Getzler order} of $Q$. Moreover, extending the 
Definition~\ref{def:Greiner.homogeneous-PsiDO} to homogeneous symbols with coefficients in $(\End \C^{p})\otimes \Lambda(n)$,
we define the \emph{model operator} of $Q$ as the element of $\Psi_{\op{v}}^{*}(\R^{n}\times\R, \C^{p})\otimes \Lambda(n)$ 
given by
\begin{equation}
    Q_{(m)}:=q_{(m)}(x,D_{x},D_{t}).
\end{equation}

In the sequel we will write $\op{O_{G}}(m)$ to denote a Volterra \psido\ of Getzler order~$\leq m$ and we will write $x\op{O_{G}}(m)$ to denote an 
element of  $\Psi_{\op{v}}^{*}(\R^{n}\times\R, \C^{p})\otimes \Cl(n)$ of the form $x^{i}Q_{i}$ with $Q_{i}$ of Getzler order~$\leq m$. 
% Notice that any Volterra \psido\ of the form 
% $x\op{O_{G}}(m)$ has Getzler order~$\leq m-1$. 

Let $A=A_{i}dx^i$ be the connection one-form on $\C^p$. Then  by~(\ref{eq:AS-asymptotic-geometric-data}) the covariant 
derivative $ \nabla_{i}= \partial_{i}+ \frac14 \omega_{ikl}(x)c(e^k)c(e^l) +A_{i}$  on $\Cl(n)\otimes \C^p$
    has Getzler order 1 and model operator,  
    \begin{equation}
        \nabla_{i (1)} =\partial_{i}-\frac14 R_{ij}^{M}(0)x^j , \qquad 
        R_{ij}^{M}(0)=R_{ijkl}^{M}(0)dx^k \wedge dx^l. 
        \label{eq:AS.model-spin-connection}
    \end{equation}
This implies that $\clD=c(dx^{i}) \nabla_{i}$ has Getzler order $2$ and can be expanded as 
\begin{equation}
    \clD= c(\clD_{(2)}) + x\op{O}_{G}(2), \qquad \clD_{(2)}=\varepsilon(dx^{i})  \nabla_{i (1)}, 
     \label{eq:BFA.sD-model}
\end{equation}
where $\varepsilon(dx^{i})$ denotes the exterior multiplication by $dx^{i}$.    

The interest to introduce Getzler orders stems from the two lemmas below. 
\begin{lemma}[{\cite[Lem.~3]{Po:CMP1}}]\label{lem:BFA.approximation-asymptotic-kernel}
Let $Q\in \Psi_{\op{v}}^{*}(\R^{n}\times\R, \C^{p})\otimes \Cl(n)$ have Getzler order $m$ and model operator 
     $Q_{(m)}$. 
     
   - If $m$ is odd, then as  $t\rightarrow 0^{+}$  we have  
      \begin{equation}
          \sigma[K_{Q}(0,0,t)]^{(n)}= t^{-(\frac{m}2+1)} 
         K_{Q_{(m)}}(0,0,1)^{(n)} + \op{O}(t^{-\frac{m}2}).
         \label{eq:BFA.approximation-odd}
      \end{equation}
      - If $m$ is even, then $K_{Q_{(m)}}(0,0,t)^{(n)}=0$ and as $t\rightarrow 0^{+}$  we have  
      \begin{equation}
          \sigma[K_{Q}(0,0,t)]^{(n)}=\op{O}(t^{-\frac{m+1}2}). 
           \label{eq:BFA.approximation-even}
      \end{equation}
\end{lemma}
\begin{proof}
 Let $q(x,\xi,\tau)\sim \sum_{j \leq m'} q_{j}(x,\xi,\tau)$ be the symbol of $Q$ and let $q_{(m)}(x,\xi,\tau)$ be its principal 
 Getzler homogeneous symbol. Then by Lemma~\ref{lem:volterra.key-lemma} as $t\rightarrow 0^{+}$ we have 
$\sigma[K_{Q}(0,0,t)]^{(n)}\sim  \sum_{j \leq m'}  t^{-\frac{n+2-j}{2}} \sigma[\check{q}_{j}(0,0,1)]^{(n)}$.
%  \begin{equation}
% 
%  \end{equation}
Notice that $\sigma[q_{j}(0,\xi,\tau)]^{(n)}$ is Getzler homogeneous of degree $j+n$, so it must be zero if $j+n>m$, since 
otherwise $Q$ would have Getzler 
order~$>m$. Therefore, we get: 
\begin{multline}
    \sigma[K_{Q}(0,0,t)]^{(n)} \\
    =t^{-(\frac{m}2+1)} 
         \sigma[\check{q}_{m-n}(0,0,1)]^{(n)} +t^{-\frac{m+1}2}
         \sigma[\check{q}_{m-n-1}(0,0,1)]^{(n)} +  \op{O}(t^{-\frac{m}2}) . 
\end{multline}

On the other hand, in view of asymptotics~(\ref{eq:AS.Getzler-asymptotic})--(\ref{eq:index.asymptotic-symbol}) 
the symbol $q_{(m)}(0,\xi,\tau)^{(n)}$ is equal to
$\sum_{j+n-|\alpha|=m}  (\frac{x^{\alpha}}{\alpha!}  
     \partial_{x}^{\alpha}\sigma[q_{j}(0,\xi,\tau)]^{(n)})_{x=0}= \sigma[q_{m-n}(0,\xi,\tau)]^{(n)}$.
% It follows that we have $\sigma[\check{q}_{m-n}(0,0,1)]^{(n)}=K_{Q_{(m)}}(0,0,1)^{(n)}$. 
Thus, 
\begin{equation}
    \sigma[\check{q}_{m-n}(0,0,1)]^{(n)}=\check{q}_{(m)}(0,0,1)^{(n)}=K_{Q_{(m)}}(0,0,1)^{(n)}.
%     \label{eq:¥}
\end{equation}
% 
% so that we get
%       \begin{equation}
%           \sigma[K_{Q}(0,0,t)]^{(n)}= t^{\frac{-m}2-1} 
%          K_{Q_{(m)}}(0,0,1)^{(n)} + \op{O}(t^{\frac{-m}2}). 
%            \label{eq:BFA.approximation2}
%       \end{equation}

Finally, as pointed out in the proof  of Lemma~\ref{lem:volterra.key-lemma} we have $\check{q}_{j}(0,0,1)=0$ when $j$ is odd. As $n$ is 
odd we see that if $m$ odd we have $ \sigma[\check{q}_{m-n-1}(0,0,1)]^{(n)}=0$ and so~(\ref{eq:BFA.approximation-odd}) holds, 
while if $m$ even then
$ K_{Q_{(m)}}(0,0,1)^{(n)}=\sigma[\check{q}_{m-n}(0,0,1)]^{(n)}=0$ and we obtain~(\ref{eq:BFA.approximation-even}). 
%       This proves the lemma when $m$ is odd. In case $m$ is even the Volterra symbol $q_{m-n}(0,\xi,\tau)$ is homogeneous of odd degree so, as alluded to 
% in the proof of , we have $\check{q}_{m-n}(0,0,1)=0$. Hence the leading term 
% in~(\ref{eq:BFA.approximation2}) vanishes and the asymptotics holds in the even case. 
\end{proof}

\begin{lemma}\label{lem:index.top-total-order-symbol-composition}
 For  $j=1,2$ let $Q_{j}\in \Psi^{*}_{\op{v}}(\R^{n}\times\R,\C^{p})\otimes \Cl(n)$ have Getzler order $m_{j}$ 
and model operator $Q_{(m_{j})}$, and assume that $Q_{1}$ or $Q_{2}$ is properly supported. Then in 
$\Psi^{*}_{\op{v}}(\R^{n}\times\R,\C^{p})\otimes \Cl(n)$ we have
\begin{gather}
    Q_{1}Q_{2}= c[Q_{(m_{1})} Q_{(m_{2})}] +\op{O}_{G}(m_{1}+m_{2}-1), 
     \label{eq:approximation-product-Volterra-PsiDO's1} \\
   c(Q_{(m_{1})})c( Q_{(m_{2})})=   c[Q_{(m_{1})} Q_{(m_{2})}] +\op{O}_{G}(m_{1}+m_{2}-2).
    \label{eq:approximation-product-Volterra-PsiDO's2}
\end{gather}
\end{lemma}
\begin{proof}
    The equality~(\ref{eq:approximation-product-Volterra-PsiDO's1}) is the content of Lemma 4 of~\cite{Po:CMP1}, 
    so we need only to prove~(\ref{eq:approximation-product-Volterra-PsiDO's2}). For for a subset 
    $I=\{i_{1},\ldots,i_{k}\}$ of  $\{1,\ldots,n\}$ with $i_{1}<\ldots<i_{k}$ we let $dx^{I}=dx^{1}\wedge \ldots \wedge dx^{n}$. Since the forms $dx^{I}$ 
    gives rise to a 
    linear basis of $\Lambda(n)$  we can write $Q_{(m_{j})}$, $j=1,2$, in the form $Q_{(m_{j})}=\sum Q_{j,I}dx^{I}$ with 
    $Q_{j,I}$ in $\Psi^{*}_{\op{v}}(\R^{n}\times\R,\C^{p})$ of Getzler order $m_{j}-|I|$. Then we have 
    \begin{equation}
        c(Q_{(m_{1})})c( Q_{(m_{2})})=\sum_{I_{1},I_{2}}Q_{1,I_{1}}Q_{2,I_{2}}c(dx^{I_{1}})c(dx^{I_{2}}).
%         \label{eq:¥}
    \end{equation}
    Thanks to~(\ref{eq:getzler.clifford-multiplication}) we have $c(dx^{I_{1}})c(dx^{I_{2}})=c(dx^{I_{1}})\wedge dx^{I_{2}})+
    \op{O}_{G}(|I_{1}|+|I_{2}|-2)$, so we obtain
    \begin{multline}
         c(Q_{(m_{1})})c( Q_{(m_{2})})=
         \sum_{I_{1},I_{2}}Q_{1,I_{1}}Q_{2,I_{2}} (c(dx^{I_{1}})\wedge dx^{I_{2}})+\op{O}_{G}(|I_{1}|+|I_{2}|-2)\\
         =c[Q_{(m_{1})} Q_{(m_{2})}] +\op{O}_{G}(m_{1}+m_{2}-2).
%          \label{eq:¥}
    \end{multline}
    The proof is thus achieved.
\end{proof}

Now, by the Lichnerowicz's formula~(\ref{eq:BFA.Lichnerowicz}) we have 
\begin{equation}
    \clD^{2}= -g^{ij}(\nabla_{i}\nabla_{j} -\Gamma_{ij}^k 
     \nabla_{k}) 
        + \frac12 c(e^i)c(e^j)F(e_{i},e_{j})  + \frac{\kappa}4,
     \label{eq:AS.lichnerowicz-bis}
\end{equation}
where the $\Gamma_{ij}^k$'s are the Christoffel symbols of the metric $g$. Observe that:\smallskip 

- Thanks to~(\ref{eq:AS-asymptotic-geometric-data}) we have $g_{ij}(x)=\delta_{ij}+\op{O}(|x|^{2})=\delta_{ij}+\op{O}_{G}(-2)$;\smallskip 

- Using (\ref{eq:approximation-product-Volterra-PsiDO's1}) and~(\ref{eq:AS.model-spin-connection}), as well as the fact that 
$\Gamma_{ij}^k(x)=\op{O}(|x|)=\op{O}_{G}(-1)$, we get $\nabla_{i}\nabla_{j}=c(\nabla_{i(1)}\nabla_{j(1)})+x\op{O}_{G}(2)+\op{O}_{G}(0)$;\smallskip

- We have $F^{\cE}(e^{i},e^{j})=F^\cE(\partial_{i},\partial_{j})(0)+\op{O}(|x|)=F^\cE(\partial_{i},\partial_{j})(0)+x\op{O}_{G}(0)$;\smallskip

- By~(\ref{eq:getzler.clifford-multiplication}) we have $c(e^i)c(e^j)= c(e^{i}\wedge e^{j})+\op{O}_{G}(0)=c(dx^{i}\wedge 
dx^{j})+x\op{O}_{G}(2)+\op{O}_{G}(0)$.\smallskip   

Moreover, as $x^{i}$ commutes with $x^{j}$ and $c(dx^{j})$ and commutes with Volterra \psidos\ of 
order $m$ modulo those of order~$\leq m-1$, we see that the commutation with $x^{i}$  preserves 
the Getzler order. In particular, if $Q\in  \Psi_{\op{v}}^{*}(\R^{n}\times\R, \C^{p})\otimes \Cl(n)$ has
Getzler order~$\leq m$, then it follows from the equality $Qx^{i}=x^{i}Q+[Q,x^{i}]$ that $Qx^{i}$ is of the form 
$x\op{O}_{G}(m)+\op{O}_{G}(m)$. 

Bearing all this in mind we obtain: 
\begin{multline}
\clD^{2}= (\delta_{ij}+ \op{O}_{G}(-2))(c(\nabla_{i(1)}\nabla_{j(1)}) +x\op{O}_{G}(2)+\op{O}_{G}(0))\\ 
+  \frac12(F^\cE(\partial_{k},\partial_{l})(0) +x\op{O}_{G}(0))(c(dx^{i}\wedge dx^{j}) +x\op{O}_{G}(2)  +\op{O}_{G}(0)) +\op{O}_{G}(0),\\ 
= c(\clD_{(2)}^{2}) +x\op{O}_{G}(2)+\op{O}_{G}(0), \qquad \clD_{(2)}^{2}=H_{R} +F^\cE(0),
\label{eq:BFA.model-sD2}
\end{multline}
where $H_{R}=- \sum_{i=1}^n (\partial_{i}-\frac14 R_{ij}^{M}(0)x^j)^{2}$ and 
$F^\cE(0)=\frac12 F^\cE(\partial_{k},\partial_{l})(0)dx^k \wedge dx^l$. In particular, we see that $\clD^{2}$ has Getzler order 2.  

\begin{lemma}\label{lem:getzler-model-asymptotic}
Any Volterra \psido\ parametrix $Q$ for $\clD^{2}+\partial_{t}$ is of the form
\begin{equation}
    \clQ=c(\clQ_{(-2)}) +x\op{O}_{G}(-2)+\op{O}_{G}(-4), \qquad \clQ_{(-2)}=(H_{R}+F^\cE(0)+\partial_{t})^{-1}.
     \label{eq:BFA.model-Q}
\end{equation}
In particular $Q$ has Getzler order $-2$.
\end{lemma}
\begin{proof}
    Thanks to~\cite[Lem.~5]{Po:CMP1} we know that $\clQ$ has Getzler order $-2$ and model operator 
    $Q_{(-2)}=(H_{R}+F^\cE(0)+\partial_{t})^{-1}$. In order to 
    get~(\ref{eq:BFA.model-Q}) notice that since
    $\clQ$ is an inverse for $\clD^{2}+\partial_{t}$ modulo  $\Psi^{-\infty}_{\op{v}}(\R^{n}\times 
    \R, \C^{p})\otimes \Cl(n)$, hence modulo operators of Getzler orders $-\infty$, we have 
    \begin{equation}
        \clQ(\clD^{2}+\partial_{t})c(\clQ_{(-2)})=(1+\op{O}_{G}(-\infty))c(Q_{(-2)})= c(Q_{(-2)}) +\op{O}_{G}(-\infty).
         \label{eq:BFA.QQ-2.1}
    \end{equation}
    
    On the other hand, from~(\ref{eq:approximation-product-Volterra-PsiDO's2}) and~(\ref{eq:BFA.model-sD2}) we get 
    \begin{multline}
        c(Q_{(-2)})(\clD^{2}+\partial_{t})=c((\clD^{2}_{(2)}+\partial_{t})^{-1})(c(\clD_{(2)}^{2}+\partial_{t}) +x\op{O}_{G}(2)+\op{O}_{G}(0))\\
        = 1+ x\op{O}_{G}(0)+\op{O}_{G}(-2).
%         \label{eq:}
    \end{multline}
    Therefore, we also obtain
%     Moreover, by we have 
%     $c(\clD^{2}_{(2)}+\partial_{t})c((\clD^{2}_{(2)}+\partial_{t})^{-1})=1+\op{O}_{G}(-2)$. Therefore, 
%     using~(\ref{eq:BFA.model-sD2}) we see that $ \clQ(\clD^{2}+\partial_{t})c(\clQ_{(-2)})$ is also equal to
    \begin{equation}
      \clQ(c(\clD^{2}_{(2)}+\partial_{t})+ 
       \op{O}_{G}(0))c((\clD^{2}_{(2)}+\partial_{t})^{-1})
        =\clQ+ x\op{O}_{G}(-2)+\op{O}_{G}(-4).
        \label{eq:BFA.QQ-2.2}
    \end{equation}
    Comparing this to~(\ref{eq:BFA.QQ-2.1}) gives $ \clQ=c(\clQ_{(-2)}) +x\op{O}_{G}(-2)+\op{O}_{G}(-4)$ as desired.
\end{proof}

Let $\clQ$ be a Volterra \psido\ parametrix for $\clD^{2}+\partial_{t}$. Combining~(\ref{eq:BFA.sD-model}), 
(\ref{eq:approximation-product-Volterra-PsiDO's2})  
and~(\ref{eq:BFA.model-Q}) we see that $\sD Q$ is equal to
\begin{multline}
    (c(\clD_{(2)})+x\op{O}_{G}(2))(c(\clQ_{(-2)})+x\op{O}_{G}(-2)+\op{O}_{G}(-4))\\ =c(\clD_{(2)}\clQ_{(-2)}) 
    +x^{j}\clQ_{j}+\clR,
\end{multline}
with $\clQ_{j}$ of Getzler order~$\leq 0$  and $\clR$ of Getzler order~$\leq-2$. 
Since $K_{x^{j}\clQ_{j}}(0,0,t)=(x^{j}K_{\clQ_{j}})(0,0,t)=0$, we obtain
\begin{equation}
      \sigma[K_{\sD \clQ}(0,0,t)]^{(n)}= K_{\clD_{(2)}\clQ_{(-2)}}(0,0,t)^{(n)}+\sigma[K_{\clR}(0,0,t)]^{(n)}.
%     \label{eq:¥}
\end{equation}

On the other hand, as $\clD_{(2)}\clQ_{(-2)}$ is Getzler homogeneous of even degree and $\clR$ has Getzler order~$\leq -2$,
Lemma~\ref{lem:BFA.approximation-asymptotic-kernel} shows that
$K_{\clD_{(2)}\clQ_{(-2)}}(0,0,t)^{(n)}=0$ and 
$\sigma[K_{\clR}(0,0,t)]^{(n)}=\op{O}(\sqrt{t})$ as $t \rightarrow 0^{+}$. It follows that as $t \rightarrow 0^{+}$ we have 
\begin{equation}
    \sigma[K_{\clQ}(0,0,t)]^{(n)}= K_{\clD_{(2)}\clQ_{(-2)}}(0,0,t)+\sigma[K_{\clR}(0,0,t)]^{(n)}=\op{O}(\sqrt{t}). 
\end{equation}
Combining this with~(\ref{eq:BFA.trace-ht-sigma}) then gives 
\begin{equation}
    \Tr_{\cE} h_{t}(x_{0},x_{0})=\op{O}(\sqrt{t}) \qquad \text{as $t \rightarrow 0^{+}$.}
\end{equation}
This shows that in the asymptotics~(\ref{eq:Greiner.heat-kernel-asymptotics}) for $ \Tr_{\cE} h_{t}(x,x)$ all the coefficients of $t^{j}$ with 
$j<\frac{1}{2}$ are 
in fact zero, that is, we recover Bismut-Freed's asymptotics~(\ref{eq:BFA.main-asymptotics}). 
% 
% can now be completed by comparing the above asymptotics  
% with the heat kernel asymptotics.  

\end{document}